\documentclass{amsart}

\usepackage{amsmath}
\usepackage{latexsym,amsfonts,amssymb,mathrsfs}
\usepackage{mathdots}
\usepackage{color}
\usepackage [all,2cell,color]{xy}
\usepackage{tikz-cd}

\UseAllTwocells
\usepackage[
textwidth=3cm,
textsize=small,
colorinlistoftodos,
]
{todonotes}

\usepackage{hyperref}

\usetikzlibrary{decorations.pathreplacing}%need this to draw nice curly brackets
\usetikzlibrary{arrows, decorations.markings} %need this for nice arrows in tikz pictures, markings in particular for double arrows
\pgfarrowsdeclarecombine{twotriang}{twotriang}{stealth}{stealth}%
{stealth}{stealth}
\tikzset{mono/.style={>-stealth}} %this is the style for the horizontal arrows
\tikzset{epi/.style={-twotriang}} %this is the style for the vertical arrows
\tikzset{twoarrowlonger/.style={double,double distance=1.5pt,
shorten <=5pt,shorten >=6pt,
decoration={markings,mark=at position -4pt with {\arrow[scale=1.75]{>}}},
preaction={decorate}}} 

\tikzset{mapstikz/.style={-stealth, 
decoration={markings,mark=at position 0pt with {\arrow[scale=0.5]{|}}}, preaction={decorate}}}%this is the style for mapsto in tikz

\tikzset{
    dot/.style={circle,draw,fill,inner sep=1pt}
}

\usepackage[capitalise]{cleveref}

\newcommand{\newrefformat}[2]{}

\theoremstyle{plain}   % This is the default, anyway

\numberwithin{equation}{section}

 \newtheorem{thm}[equation]{Theorem}

\newtheorem*{thm*}{Theorem}
\newtheorem*{cor*}{Corollary}
\newtheorem*{lem*}{Lemma}
\newtheorem*{prop*}{Proposition}
\newtheorem*{claim*}{Claim}

\theoremstyle{definition}
\newtheorem{defn}[equation]{Definition}
\newtheorem*{defn*}{Definition}

\theoremstyle{remark}
\newtheorem{rmk}[equation]{Remark}
\newtheorem{ex}[equation]{Example}

\newtheorem*{rmk*}{Remark}
\newtheorem*{ex*}{Example}
\newtheorem*{exs*}{Examples}

\newcommand{\cA}{\mathcal{A}}
\newcommand{\cB}{\mathcal{B}}
\newcommand{\cC}{\mathcal{C}}

\newcommand{\cE}{\mathcal{E}}

\newcommand{\cM}{\mathcal{M}}

\newcommand{\cQ}{\mathcal{Q}}

\newcommand{\sdot}{S_{\bullet}}%macro for Waldhausen Sdot construction

\newcommand{\htimes}[1]{\underset{#1}{\overset{h}{\times}}}

\newcommand{\incls}[2]{\alpha^{#1}_{#2}}
\newcommand{\inclt}[2]{\widetilde{\alpha}^{#1}_{#2}}

\newcommand{\funsp}{P}

\DeclareMathOperator{\esd}{esd} %edgewise subdivision

\DeclareMathOperator{\Fun}{Fun}

\begin{document}
\title{The edgewise subdivision criterion for $2$-Segal objects}

\author[J.~E.~Bergner]{Julia E.~Bergner}
\address{Department of Mathematics, University of Virginia, Charlottesville, VA 22904, USA}
\email{jeb2md@virginia.edu}
\author[A.~M.~Osorno]{Ang\'{e}lica M. Osorno}
\address{Department of Mathematics, Reed College, Portland, OR 97202, USA}
\email{aosorno@reed.edu}
\author[V.~Ozornova]{Viktoriya Ozornova}
\address{Fakult\"at f\"ur Mathematik, Ruhr-Universit\"at Bochum, D-44780 Bochum, Germany}
\email{viktoriya.ozornova@rub.de}

\author[M.~Rovelli]{Martina Rovelli}
\address{Department of Mathematics,
Johns Hopkins University,
Baltimore, MD 21218, USA}
\email{mrovelli@math.jhu.edu}

\author[C.~Scheimbauer]{Claudia I.~Scheimbauer}
\address{Mathematical Institute, University of Oxford, Oxford OX2 6GG, UK}
\email{scheimbauer@maths.ox.ac.uk}

\thanks{The first-named author was partially supported by NSF CAREER award DMS-1352298. The second-named author was partially supported by a grant from the Simons Foundation (\#359449, Ang\'elica Osorno) and NSF grant DMS-1709302.
The fourth-named author and fifth-named authors were partially funded by the Swiss National Science Foundation, grants P2ELP2\textunderscore172086 and P300P2\textunderscore 164652, respectively.}

\begin{abstract}
We show that the edgewise subdivision of a $2$-Segal object is always a Segal object, and furthermore that this property characterizes $2$-Segal objects.
\end{abstract}

\maketitle

\section{Introduction}

The edgewise subdivision of a simplicial space is a construction which leaves the geometric realization unchanged but has the effect of decomposing the simplicial space into more simplices.  It first appeared in the literature in work of Segal \cite{segal}, although he attributes it to Quillen.   Waldhausen in turn used this construction to prove the equivalence of the $S_\bullet$-construction and the $Q$-construction in algebraic $K$-theory \cite{waldhausen}.

In this note, we consider the effect of applying this construction to the $2$-Segal spaces of Dyckerhoff and Kapranov \cite{DK} (closely related to the decomposition spaces of G\'alvez-Carrillo, Kock, and Tonks \cite{GalvezKockTonks}), which model homotopical categories with associative multivalued composition.  A key source of examples of such structures is the output of Waldhausen's $\sdot$-construction when applied to an exact category, as shown in \cite{DK} and \cite{GalvezKockTonks}.  In \cite{BOORS} and \cite{BOORS2} we show that any 2-Segal space which satisfies a unitality condition arises from such a construction for a suitably general input.
 
In this paper we work in the more general context of $2$-Segal objects in any combinatorial model category.  Our main result, which appears as \cref{mainresult}, characterizes 2-Segal objects in terms of their edgewise subdivision.

\begin{thm*}
 Let $X$be a simplicial object in a combinatorial model category $\cM$.  Then $X$ is a $2$-Segal object if and only if its edgewise subdivision $\esd(X)$ is a Segal object.
\end{thm*}

Although we prove the result as stated, the nature of our arguments is purely combinatorial and the statement should also hold in more general settings, such as for $2$-Segal objects in any $(\infty,1)$-category with finite limits.

This criterion will be used in forthcoming work of G\'alvez-Carrillo, Kock and Tonks \cite{gktnew}.

Following a background section and the proof of the main theorem, we conclude the paper with additional examples. In particular, we prove that the $\sdot$-construction of the Waldhausen category of retractive spaces is not a $2$-Segal space.

\subsection*{Acknowledgements}

We would like to thank the organizers of the Women in Topology II Workshop and the Banff International Research Station for providing a wonderful opportunity for collaborative research, during which this project was started.

\section{The main theorem}

In this section, we recall the necessary background to state our main theorem: the definition of the edgewise subdivision of a simplicial object and an overview of $2$-Segal objects.

Let $\Delta$ be the simplicial indexing category, whose objects are non-empty finite ordered sets and whose morphisms are order-preserving maps. Recall that the \emph{join} $\cA \star \cB$ of two categories $\cA$ and $\cB$ is obtained by taking the disjoint union of the two categories, then adjoining a unique additional arrow from every object of $\cA$ to every object in $\cB$ \cite[\S 1.2.8]{htt}. 

We follow the convention of \cite{waldhausen} and denote by $\epsilon \colon \Delta \to \Delta$ the functor determined by the assignment  $[n]\mapsto[n]^\text{op}\star [n]\cong[2n+1]$.
The image of $[n]$ can be depicted as the ordered set
\[n' < (n-1)' < \cdots < 1' < 0' < 0 < 1 < \cdots < n-1 < n.\] 

\begin{defn} \label{defn: edge sub}
Let $\cC$ be any category. The \emph{edgewise subdivision} is the functor
$$\esd\colon\Fun(\Delta^{\text{op}},\cC)\to\Fun(\Delta^{\text{op}},\cC)$$
that assigns to a simplicial object $X$ the simplicial object $\esd(X)=X \circ \epsilon$, and is therefore given in component $n$ by $\esd(X)_n =X_{2n+1}$.
\end{defn}

We warn the reader of different conventions in the literature; we chose to follow \cite[\S 1.9]{waldhausen} and \cite[Example 2.5]{barwickq}, but other references, such as \cite[Appendix 1]{segal} and \cite[\S 10]{DK}, use the opposite convention of $[n] \star [n]^{\text{op}}$ instead.
See \cite{Velcheva} for further discussions of variants of subdivisions, some of which are fundamentally different, for example the one used initially in \cite{GraysonExterior}. 
The terminology is justified by the following example.

\begin{ex}
The edgewise subdivision $\esd(\Delta[k])$ of a standard simplex $\Delta[k]$ is indeed a subdivision of $\Delta[k]$ into $2^k$ non-degenerate $k$-simplices, as was shown in Segal's work \cite{segal}. For example, the subdivision of $\Delta[2]$ given by $\esd(\Delta[2])$ is depicted as
$$
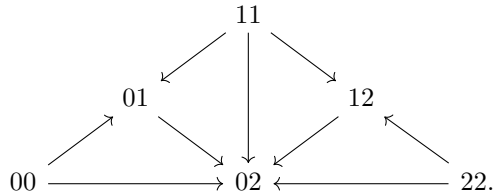
\begin{tikzcd}
&&11\\
&01\arrow[from=ur]\arrow[from=dl] && 12 \arrow[from=ul]\arrow[from=dr]\\
00 && 02 \arrow[from=ll]\arrow[from=ul]\arrow[from=uu]\arrow[from=ur]\arrow[from=rr] && 22.
\end{tikzcd}
$$
\end{ex}

\begin{ex}
The edgewise subdivision of the nerve of a category $\mathcal{A}$ is isomorphic to the nerve of an associated familiar category.  Recall that, given a category $\mathcal{A}$, its \emph{twisted arrow category} $\mathrm{Tw}(\cA)$ is the category whose objects are the morphisms $f\colon a\to b$ in $\cA$, and whose morphisms from $f$ to $g$ are given by commutative squares
$$\begin{tikzcd}
a  \arrow[swap]{d}{f} & c \arrow{l} \arrow{d}{g}\\
b \arrow{r} & d. 
\end{tikzcd}$$
Then there is an isomorphism of simplicial sets
$$\esd(N\cA)\cong N \mathrm{Tw}(\cA),$$
as explained by Barwick in \cite[\S 2.6]{barwickq}, specializing some of the ideas and constructions from the proof of Quillen's Theorem A \cite{QuillenK}.

This example serves as inspiration for the definition of twisted arrow quasi-category (see \cite[\textsection5.2.1]{LurieHA} and \cite[\textsection2.6]{barwickq}, with some differing conventions), which is defined precisely as the simplicial set $\esd(\cQ)$ for any quasi-category $\cQ$.
\end{ex}

In particular, the previous example shows that the edgewise subdivision  of the nerve of a category is always given by the nerve of a category. An example exhibiting a similar phenomenon, due to Segal, is that of partial topological monoids, which we now recall.

\begin{ex}  \label{ex: partial monoids}
In \cite{segal}, Segal introduced partial topological monoids and their classifying spaces to relate configuration spaces and iterated loop spaces.  Recall that a \emph{partial topological monoid} is a space $M$ together with a subspace $M_2\subset M\times M$ and a map $M_2\to M$, written as $(m_1,m_2)\mapsto m_1\cdot m_2$, which we think of as a multiplication when it is defined. This partially defined multiplication is assumed to be unital and associative whenever both sides of the usual associativity condition are defined. 

Segal defines the nerve of $M$ to be the simplicial space $BM=M_\bullet$ where $M_0=*$ and
$$M_n= \{(m_1,\ldots, m_n) \in M^{\times n}\colon \text{ for } 1\leq i< n, (m_1\cdots m_i, m_{i+1})\in M_2\},$$
with faces and degeneracies given by partial multiplication and insertions of the unit, respectively.  Moreover, he defines a topological category $\mathscr{C}(M)$ whose objects are the elements of $M$ and whose morphisms $m \rightarrow m'$ are given by pairs of elements $m_1,m_2\in M$ such that $m_1\cdot m \cdot m_2=m'$. In particular, every morphism can be identified with a triple $(m_1,m,m_2)\in M_3$. Using this formulation, we can express composition of morphisms $(m_1,m, m_2)$ and $(n_1, m_1\cdot m \cdot m_2, n_2)$ by $(n_1 \cdot m_1, m, m_2\cdot n_2)$. Note that triple is a well-defined element in $M_3$ since $(n_1 \cdot m_1)\cdot m \cdot(m_2 \cdot n_2) = n_1\cdot(m_1\cdot m \cdot m_2)\cdot n_2 \in M$. Associativity of composition follows from associativity of the partial multiplication in the monoid. 

One can verify that these two constructions are related precisely by edgewise subdivision, via an isomorphism 
\begin{equation}\label{eq: partial monoid}
\esd(BM)\cong N\mathscr{C}(M).
\end{equation}
Again, we see that in this case the output of the edgewise subdivision functor is the nerve of a category. Note after that taking geometric realization on both sides this is precisely the statement in \cite[Proposition 2.5]{segal}. 
\end{ex} 

Our main result is that this feature of both examples, that the edgewise subdivision can be described as a nerve of a (possibly topological) category, is not a coincidence.  To state this result, however, we need to describe the appropriate structure for the input of the edgewise subdivision functor, that of 2-Segal objects.  We now review this structure, referring the reader to \cite{DK} for more detailed constructions and proofs.

Segal and $2$-Segal objects are often considered in the category of simplicial sets with its original model structure, due to Quillen, (Theorem 3 in \cite[\S II.3]{QuillenHA}).
However, since our main result holds in more generality, we recall the definitions of Segal and $2$-Segal objects in a combinatorial model category $\cM$ as e.g.\ in \cite[A.2.6]{htt}. We choose to work in this setting to simplify technical arguments, but we could instead work in a more general setting, such as in an $(\infty,1)$-category $\cM$ with finite limits, see \cref{comment}.
 
The definitions which we give here are not those from the original sources, namely \cite[\textsection4]{rezk} and \cite[Definition 2.3.1]{DK}, respectively, but they are equivalent and are better suited for the purposes of this paper.

\begin{defn}
A \emph{Segal object} in $\cM$ is a simplicial object in $\cM$ such that, for every $m\geq j \geq 1$, the Segal map
\[\beta^m_j\colon X_m \rightarrow X_j \htimes{X_0} X_{m-j}, \]
induced by the maps
$\incls{j}{m} \colon [j] \to [m]$, given by $i \mapsto i$ for all $0 \leq i \leq j$, and $\inclt{j}{m} \colon [m-j] \to [m]$, given by $i \mapsto i+j$, for all $0 \leq i \leq m-j$, is a weak equivalence.
\end{defn}
Graphically, this map corresponds to the decomposition of the interval $[0,m]$ with $m$ segments into the intervals $[0,j]$ and $[j,m]$,
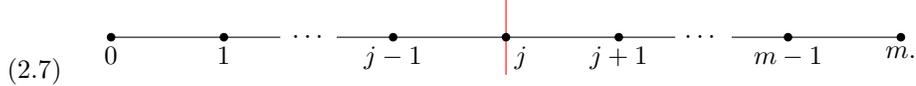
\begin{equation}\label{pic 1-Segal}
\begin{tikzpicture}[xscale=1.5]
\draw (0,0) node[dot](0) {} node[anchor=north] {$0$} -- (1,0)node[dot](1) {} node[anchor=north] {$1$} -- (1.5, 0);

\begin{scope}[xscale=-1, xshift=-7cm]
\draw (0,0) node[dot](0) {} node[anchor=north] {$m.$} -- (1,0)node[dot](1) {} node[anchor=north] {$m-1$} -- (1.5,0);
\end{scope}

\draw[color=red] (3.5, 0.5) -- (3.5, -0.5);
\draw (2,0) node[anchor=east] {$\cdots$} --  (2.5,0) node[dot] (j-1) {} node[anchor=north]{$j-1$}
-- (3.5, 0) node[dot] (j) {} node[anchor=north, xshift=0.2cm] {$j$}
-- (4.5, 0) node[dot] (j+1) {} node[anchor=north]{$j+1$} -- (5,0) node[anchor=west] {$\cdots$};
\end{tikzpicture}
\end{equation}

In \cite{DK}, Dyckerhoff and Kapranov generalized this notion to that of $2$-Segal objects, for which the maps above are generalized to ones induced by polygonal decompositions of polygons. Given any $(n+1)$-gon and any decomposition $\mathcal T$ of that polygon by a single diagonal between vertices indexed by $i$ and $j$, as depicted in
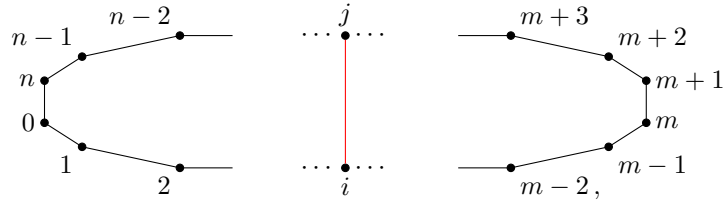
\begin{center}
\begin{tikzpicture}[yscale=0.4]
\draw (-1.5, -2.2) -- (-2.2, -2.2) node[dot] (2) {} node[anchor=north east] {$2$} -- (-3.5, -1.5) node[dot] (1) {} node[anchor=north east] {$1$} -- (-4, -0.7) node[dot] (0) {} node[anchor=east] {$0$} -- (-4, 0.7) node[dot] (0') {} node[anchor=east] {$n$} -- (-3.5, 1.5) node[dot] (1') {} node[anchor=south east] {$n-1$} -- (-2.2, 2.2) node[dot] (2') {} node[anchor=south east] {$n-2$} -- (-1.5, 2.2) ;

\draw (0, -2.2) node[dot] (j) {} node[anchor=north] {$i$}
(0, 2.2) node[dot] (j') {} node[anchor=south] {$j$};
\draw (0, -2.2) node[anchor=east] {$\cdots$}
(0, -2.2) node[anchor=west] {$\cdots$}
(0, 2.2) node[anchor=east] {$\cdots$}
(0, 2.2) node[anchor=west] {$\cdots$};

\draw (1.5, -2.2) -- (2.2, -2.2) node[dot] (m-2) {} node[anchor=north west] {$m-2$\,,} -- (3.5, -1.5) node[dot] (m-1) {} node[anchor=north west] {$m-1$} -- (4, -0.7) node[dot] (m) {} node[anchor=west] {$m$} -- (4, 0.7) node[dot] (n') {} node[anchor=west] {$m+1$} -- (3.5, 1.5) node[dot] (m-1') {} node[anchor=south west] {$m+2$} -- (2.2, 2.2) node[dot] (m-2') {} node[anchor=south west] {$m+3$} -- (1.5, 2.2) ;

\draw[color=red] (j) -- (j');

\end{tikzpicture}
\end{center}
there is a natural inclusion of simplicial sets
$$\Delta[1,\ldots, i,j,\ldots, n] \underset{\Delta[i,j]}{\amalg} \Delta[i,\ldots, j]\hookrightarrow\Delta[n].$$

We make use of the following definition of a $2$-Segal space, which is proven to be equivalent to the original one in \cite{DK}.

\begin{defn} 
A \emph{$2$-Segal object} in $\cM$ is a simplicial object such that, for every $n \geq 3$ and $0 \leq i < j \leq n$, the induced map
\[\gamma_{i,j}^n: X_n \rightarrow X_{\{0, \ldots, i, j, \ldots,n\}} \htimes{X_{\{i,j\}}} X_{\{i, \ldots, j\}} \]
is a weak equivalence.  
\end{defn}

An inductive argument can be used to show that it suffices that these maps be weak equivalences in the special cases when $i=0$ or $j=n$ (see \cite[Proposition 2.3.2]{DK}).

We can observe that the inputs to edgewise subdivision in the above examples have the common structure of $2$-Segal objects.  We can now state our main result precisely.

\begin{thm}
\label{mainresult}
 Let $X$be a simplicial object in a combinatorial model category $\cM$.  Then $X$ is $2$-Segal if and only if $\esd(X)$ is Segal.
\end{thm}

The ``if'' direction of this theorem is straightforward. However, the ``only if'' direction is more surprising. Note that it is similar to the Path Space Criterion (\cite[Theorem 6.3.2]{DK} or \cite[Proposition 4.9]{GalvezKockTonks}).  We devote the next section to the proof, after revisiting one of our examples.

\begin{ex}
The nerve of any category internal to topological spaces is a Segal object in topological spaces, so in particular $N\mathscr{C}(M)$ is. We showed in \cite[Example 2.1]{BOORS} that $BM$ is a 2-Segal topological space. Thus \eqref{eq: partial monoid} can be regarded as a special case of \cref{mainresult}.
\end{ex}

\section{Proof of the main theorem} 
We now give the proof of \cref{mainresult}.  Our strategy is to exhibit the $2$-Segal condition for $X$ in terms of the Segal condition for $\esd(X)$ and conversely, and we do so geometrically.

\begin{proof}[Proof of \cref{mainresult}]
Assume that $X$ is $2$-Segal. We need to show that $\esd(X)$ satisfies the Segal condition, namely, that for all $m\geq 2$ and $0<j<m$, the map
\[\beta^m_j\colon \esd(X)_m \longrightarrow \esd(X)_{j} \overset{h}{\underset{\esd(X)_{0}}{\times}} \esd(X)_{m-j}\]
induced by $\incls{j}{m}$ and $\inclt{j}{m}$ is a weak equivalence.

Using the definition of $\epsilon$, we can represent the maps
\begin{align*}
[2j+1]\cong\epsilon([j]) &\xrightarrow{\epsilon(\incls{j}{m})} \epsilon([m]) \cong [2m+1], \quad\text{and}\\
[2(m-j)+1]\cong\epsilon([m-j]) &\xrightarrow{\epsilon(\inclt{j}{m})} \epsilon([m]) \cong [2m+1]
\end{align*}
as 
\begin{align*}
j'<\cdots<0'<0<\cdots < j & \longmapsto  j'<\cdots<0'<0<\cdots < j \quad\text{and}\\
(m-j)'<\cdots<0'<0<\cdots < m-j &\longmapsto  m'<\cdots<j'< j<\cdots < m,
\end{align*}
respectively. Thus, $\beta^m_j$ is exactly the 2-Segal map $\gamma^{2m+1}_{m-j,m+j+1}$ for the following polygonal decomposition of the $(2m+2)$-gon, which arises as a ``cylinder'' over the one-dimensional decomposition of the interval from \eqref{pic 1-Segal}:
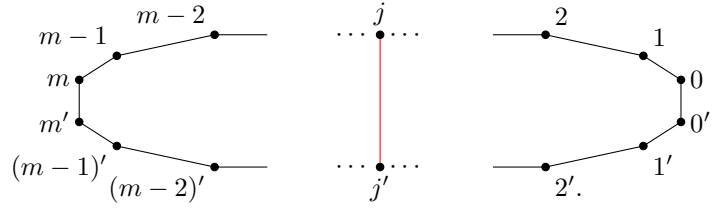
\begin{center}
\begin{tikzpicture}[yscale=-0.4, xscale=-1]
\draw (-1.5, -2.2) -- (-2.2, -2.2) node[dot] (2) {} node[anchor=south west] {$2$} -- (-3.5, -1.5) node[dot] (1) {} node[anchor=south west] {$1$} -- (-4, -0.7) node[dot] (0) {} node[anchor=west] {$0$} -- (-4, 0.7) node[dot] (0') {} node[anchor=west] {$0'$} -- (-3.5, 1.5) node[dot] (1') {} node[anchor=north west] {$1'$} -- (-2.2, 2.2) node[dot] (2') {} node[anchor=north west] {$2'.$} -- (-1.5, 2.2) ;

\draw (0, -2.2) node[dot] (j) {} node[anchor=south] {$j$}
(0, 2.2) node[dot] (j') {} node[anchor=north] {$j'$};
\draw (0, -2.2) node[anchor=west] {$\cdots$}
(0, -2.2) node[anchor=east] {$\cdots$}
(0, 2.2) node[anchor=west] {$\cdots$}
(0, 2.2) node[anchor=east] {$\cdots$};

\draw (1.5, -2.2) -- (2.2, -2.2) node[dot] (m-2) {} node[anchor=south east] {$m-2$} -- (3.5, -1.5) node[dot] (m-1) {} node[anchor=south east] {$m-1$} -- (4, -0.7) node[dot] (m) {} node[anchor=east] {$m$} -- (4, 0.7) node[dot] (n') {} node[anchor=east] {$m'$} -- (3.5, 1.5) node[dot] (m-1') {} node[anchor=north east] {$(m-1)'$} -- (2.2, 2.2) node[dot] (m-2') {} node[anchor=north east] {$(m-2)'$} -- (1.5, 2.2) ;

\draw[color=red] (j) -- (j');

\end{tikzpicture}
\end{center}
Conversely, suppose $\esd(X)$ is Segal.
Thus, we know that for subdivisions of the $(2m+2)$-gon as above
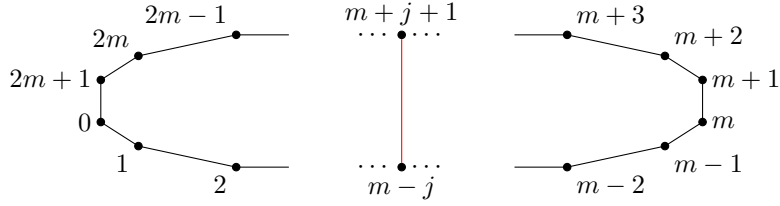
\begin{center}
\begin{tikzpicture}[yscale=0.4]
\draw (-1.5, -2.2) -- (-2.2, -2.2) node[dot] (2) {} node[anchor=north east] {$2$} -- (-3.5, -1.5) node[dot] (1) {} node[anchor=north east] {$1$} -- (-4, -0.7) node[dot] (0) {} node[anchor=east] {$0$} -- (-4, 0.7) node[dot] (0') {} node[anchor=east] {$2m+1$} -- (-3.5, 1.5) node[dot] (1') {} node[anchor=south east] {$2m$} -- (-2.2, 2.2) node[dot] (2') {} node[anchor=south east] {$2m-1$} -- (-1.5, 2.2) ;

\draw (0, -2.2) node[dot] (j) {} node[anchor=north] {$m-j$}
(0, 2.2) node[dot] (j') {} node[anchor=south] {$m+j+1$};
\draw (0, -2.2) node[anchor=east] {$\cdots$}
(0, -2.2) node[anchor=west] {$\cdots$}
(0, 2.2) node[anchor=east] {$\cdots$}
(0, 2.2) node[anchor=west] {$\cdots$};

\draw (1.5, -2.2) -- (2.2, -2.2) node[dot] (m-2) {} node[anchor=north west] {$m-2$} -- (3.5, -1.5) node[dot] (m-1) {} node[anchor=north west] {$m-1$} -- (4, -0.7) node[dot] (m) {} node[anchor=west] {$m$} -- (4, 0.7) node[dot] (n') {} node[anchor=west] {$m+1$} -- (3.5, 1.5) node[dot] (m-1') {} node[anchor=south west] {$m+2$} -- (2.2, 2.2) node[dot] (m-2') {} node[anchor=south west] {$m+3$} -- (1.5, 2.2) ;

\draw[color=red] (j) -- (j');

\end{tikzpicture}
\end{center}
the $2$-Segal map
\[\gamma^{2m+1}_{m-j,m+j+1} \colon X_{2m+1} \xrightarrow{\simeq} X_{2m-2j+1} \overset{h}{\underset{X_1}{\times}} X_{2j+1}\]
is indeed a weak equivalence.

It remains to check the $2$-Segal condition for subdivisions of an arbitrary $(n+1)$-gon into two polygons by adding a new edge that involves one of the vertices labeled with 0 or $n$. Since the other situation is symmetric, it is enough to check the cases when the new edge starts at 0. Fix $n\geq 3$ and $1<k<n$, and consider the addition of the edge between $0$ and $k$:
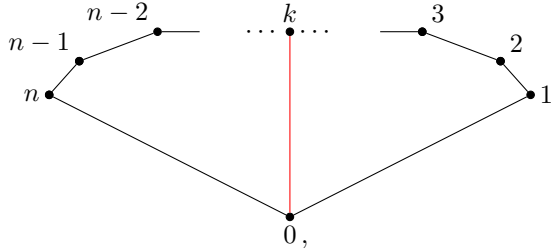
\begin{center}
\begin{tikzpicture}[yscale=0.7, scale=0.8]

\draw (4, 0.7) node[dot] (m'){};
\draw (0, -2.2) node[dot] (j) {} node[anchor=north] {} (0, 2.2) node[dot] (j') {} node[anchor=south] {};
\draw (0, 2.2) node[anchor=east] {$\cdots$} (0, 2.2) node[anchor=west] {$\cdots$};
\draw (-2.2, 2.2) node[dot] (2') {};
\draw (-4, 0.7) node[dot] (0') {};
\draw (-3.5, 1.5) node[dot] (1') {};
\draw  (3.5, 1.5) node[dot] (m-1') {}; 
\draw  (2.2, 2.2) node[dot] (m-2') {};

\draw (-1.5, 2.2) -- (2') node[anchor=south east] {$n-2$} -- (1')node[dot] (b) {} node[anchor=south east] {$n-1$} -- (0')node[dot] (c) {} node[anchor=east] {$n$}-- (j) node[anchor=north] {$\hphantom{\,,}0\,,$}-- (m') node[anchor=west] {$1$} -- (m-1') node[dot] (d) {}node[anchor=south west] {$2$} -- (m-2')node[dot] (aa) {} node[anchor=south west] {$3$} -- (1.5, 2.2);

\draw (j') node[anchor=south] {$k$};

\draw[color=red] (j) -- (j');

\end{tikzpicture}
\end{center}
i.e., we need to show that the induced map
\[\gamma^n_{0,k}\colon X_n\longrightarrow X_{n-k+1} \htimes{X_1} X_k\]
is a weak equivalence.
  
Consider the following embedding of this polygon into a larger polygon together with a decomposition which is of the above form:
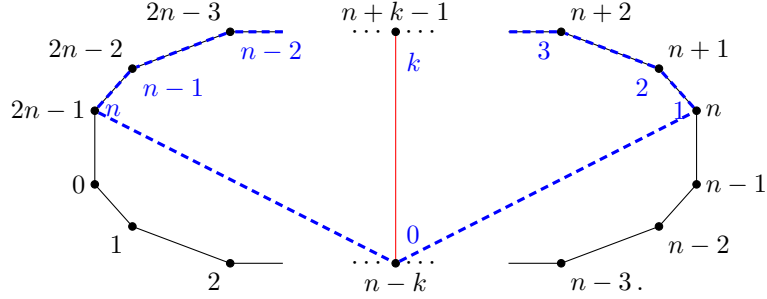
\begin{center}
\begin{tikzpicture}[yscale=0.7]
\draw (-1.5, -2.2) -- (-2.2, -2.2) node[dot] (2) {} node[anchor=north east] {$2$} -- (-3.5, -1.5) node[dot] (1) {} node[anchor=north east] {$1$} -- (-4, -0.7) node[dot] (0) {} node[anchor=east] {$0$} -- (-4, 0.7) node[dot] (0') {} node[anchor=east] {$2n-1$} -- (-3.5, 1.5) node[dot] (1') {} node[anchor=south east] {$2n-2$} -- (-2.2, 2.2) node[dot] (2') {} node[anchor=south east] {$2n-3$} -- (-1.5, 2.2) ;

\draw (0, -2.2) node[dot] (j) {} node[anchor=north] {$n-k$}
(0, 2.2) node[dot] (j') {} node[anchor=south] {$n+k-1$};
\draw (0, -2.2) node[anchor=east] {$\cdots$}
(0, -2.2) node[anchor=west] {$\cdots$}
(0, 2.2) node[anchor=east] {$\cdots$}
(0, 2.2) node[anchor=west] {$\cdots$};

\draw (1.5, -2.2) -- (2.2, -2.2) node[dot] (m-2) {} node[anchor=north west] {$n-3 \,.$} -- (3.5, -1.5) node[dot] (m-1) {} node[anchor=north west] {$n-2$} -- (4, -0.7) node[dot] (m) {} node[anchor=west] {$n-1$} -- (4, 0.7) node[dot] (m') {} node[anchor=west] {$n$} -- (3.5, 1.5) node[dot] (m-1') {} node[anchor=south west] {$n+1$} -- (2.2, 2.2) node[dot] (m-2') {} node[anchor=south west] {$n+2
$} -- (1.5, 2.2) ;

\draw[blue, very thick, densely dashed] (-1.5, 2.2) -- (2') node[anchor=north west] {$n-2$} -- (1') node[anchor=north west] {$n-1$} -- (0') node[anchor=west] {$n$}-- (j) node[label=87:0] {}-- (m') node[anchor=east] {$1$} -- (m-1') node[anchor=north east] {$2$} -- (m-2') node[anchor=north east] {$3$} -- (1.5, 2.2);

\draw[blue, very thick] (j') node[label = 273:$k$] {};

\draw[color=red] (j) -- (j');

\end{tikzpicture}
\end{center}

This decomposition determines the following commutative diagram in $\Delta$, where the bottom face is given by the decomposition of the smaller polygon, the top face is given by the decomposition of the larger polygon by the red edge (from $n-k$ to $n+k-1$), and the vertical maps describe precisely the embedding of the smaller polygon into the larger one and their decompositions. 
\begin{equation}\label{diag 1}
\begin{tikzcd}[column sep=tiny, row sep=tiny]
{[2n-1]} \arrow[rr, leftarrow] \arrow[dr,leftarrow] \arrow[dd,leftarrow, "\delta", swap, minimum width=0] &&
  {[2k-1]} 
  \\
& {[2n-2k+1]} &&
  {[1]} \arrow[dd, equal] \arrow[ul]\\
{[n]} \arrow[rr,leftarrow] \arrow[dr,leftarrow] && {[k]} \arrow[uu] \arrow[dr,leftarrow] \\
& {[n-k+1]} \arrow[rr, leftarrow] \arrow[uu, crossing over]&& {[1].}  
\ar[from=2-4, to=2-2, crossing over]
\end{tikzcd}
\end{equation}
The map $\delta$ is a composite of faces explicitly given by 
\[\delta(i)=\begin{cases}
n-k \quad &\text{if }i=0\\
i+n-1 \quad &\text{if }i>0.
\end{cases}
\]
The other vertical maps are given by similar inclusions. There is another cube with the same top and bottom faces, coming from the collapse of the non-dashed polygon edges in the inclusion:
\begin{equation}\label{diag 2}
\begin{tikzcd}[column sep=tiny, row sep=tiny]
{[2n-1]} \arrow[rr, leftarrow] \arrow[dr,leftarrow] \arrow[dd, "\sigma", swap, minimum width=0] &&
  {[2k-1]}
  \\
& {[2n-2k+1]} &&
  {[1]} \arrow[dd, equal] \arrow[ul]\\
{[n]} \arrow[rr,leftarrow] \arrow[dr,leftarrow] && {[k]} \arrow[uu, leftarrow] \arrow[dr,leftarrow] \\
& {[n-k+1]} \arrow[rr, leftarrow] \arrow[uu, crossing over, leftarrow]&& {[1].}  
\ar[from=2-4, to=2-2, crossing over]
\end{tikzcd}
\end{equation}
Here $\sigma$ is a composite of degeneracy maps given by 
\[\sigma (i)=\begin{cases}
0 \quad & \text{if }i\leq n-1\\
i-n+1 \quad  & \text{if }i\geq n.
\end{cases}\]
The other vertical maps are similarly given by composites of degeneracies.
As can been seen geometrically, $\sigma \circ \delta$ is equal to the identity map, and similarly for the corresponding composites of the other vertical maps.
Thus stacking the two commutative cubes in $\Delta$ we obtain a commutative diagram in $\cM$,
\begin{equation}\label{diag 3}
\begin{tikzcd}[column sep={4.1cm,between origins}]
X_n \arrow{r}{\gamma^n_{0,k}} \arrow{d}{\sigma^*} \arrow[bend right =80, swap]{dd}{\operatorname{id}}&  X_{n-k+1} \htimes{X_1} X_k \arrow{d} \arrow[bend left =80]{dd}{\operatorname{id}}\\
X_{2n-1} \arrow[r, "\gamma^{2n-1}_{n-k,n+k-1}", "\simeq" swap] \arrow[d, "{\delta^*}"] & X_{2n-2k+1} \htimes{X_1} X_{2k-1}\arrow[d]\\
X_n \arrow{r}{\gamma^n_{0,k}}  &  X_{n-k+1} \htimes{X_1} X_k .
\end{tikzcd}
\end{equation}
Hence $\gamma^n_{0,k}$ is a retract of the map $\gamma^{2n-1}_{n-k,n+k-1}$
which is an equivalence by assumption. It follows from the axioms of a model category that $\gamma^n_{0,k}$ is a weak equivalence as desired.
\end{proof}

\begin{rmk}\label{comment}
Note that the argument relies on the fact that stacking the diagrams \eqref{diag 1} and \eqref{diag 2} leads to a cube in $\Delta$ whose vertical maps are identities. Therefore, we can interpret \eqref{diag 3} for an $(\infty,1)$-category $\cM$ with finite limits and arrive to the same conclusion by using that retracts of equivalences are equivalences.
\end{rmk}

\section{Further examples}

In this section we give some examples, which showcase our main theorem. The first two examples, which are well-known to experts, precisely identify both sides and illustrate the statement. We conclude with a non-example arising from the category of retractive spaces.

\begin{ex}  \label{ex: exact categories}
Given an exact category $\cC=(\cC,\cM,\cE)$, as in \cite[\textsection 2]{QuillenK} or \cite[IV]{kbook}, Quillen's {\em Q-construction} is a category $Q(\cC)$
whose objects are the objects of $\cC$ and whose morphisms from an object $a$ to an object $c$ are given by equivalence classes of spans
 \begin{center}
 \begin{tikzpicture}[scale=0.5]
 \draw (0,0) node(a){$a$}; 
 \draw (2,2) node(b){$b$};
 \draw (4,0) node(c){$c,$};
 \draw[mono] (b)--node[anchor=west, yshift=0.1cm](MM){{\scriptsize $\in \cM$}}(c);
 \draw[epi] (b)--node[anchor=east, yshift=0.1cm](EE){{\scriptsize $\cE\ni$}}(a);
 \end{tikzpicture}
 \end{center}
where two such spans are identified if there is an isomorphism of middle objects both over $a$ and over $c$.
Composition of two equivalence classes of spans along a common object is given by pulling back, as displayed
 \begin{center}
 \begin{tikzpicture}[scale=0.5]
 \draw (0,0) node(a){$a$}; 
 \draw (2,2) node(b){$b$};
 \draw (4,0) node(c){$c$};
  \draw (6,2) node(d){$d$};
 \draw (8,0) node(e){$e$ \,.};
 \draw (4,4) node[inner sep=0, outer sep=0](f){$b\underset{c}{\times} d$};
 \draw[mono] (b)--(c);
 \draw[epi] (b)--(a);
 
  \draw[mono] (d)--(e);
 \draw[epi] (d)--(c);
 
  \draw[mono] (f)--(d);
 \draw[epi] (f)--(b);
 \end{tikzpicture}
 \end{center}

In \cite{waldhausen}, Waldhausen defines his $\sdot$-construction for exact categories and shows in \cite[\textsection1.9]{waldhausen} that the $Q$-construction is compatible with edgewise subdivision in that there is a levelwise weak equivalence of simplicial spaces
$$\esd(\sdot\cC)\simeq NQ(\cC),$$
where $N$ denotes what today is called Rezk's classifying diagram for categories from \cite[\textsection3.5]{rezk}, which is a Segal space. Thus, $\sdot\cC$ is a 2-Segal space, which was proven directly in \cite{DK} and was the original motivation for the definition of 2-Segal spaces.
\end{ex}

Edgewise subdivision also recovers the complete Segal spaces of spans and cospans in a stable quasi-category.

\begin{ex}  \label{ex: stable infinity categories}
Given a quasi-category with finite limits $\cQ$, there are 
several constructions for an $(\infty,1)$-category of spans in $\cQ$, such as the quasi-category of spans as constructed in \cite[\S 10]{DK}, and the complete Segal space of spans $\text{Span}(\cQ)$, of \cite[Definition 3.3]{BarwickRognes} which is also used in \cite[Theorem 1.2]{HaugsengIteratedSpans}.  One can define analogous constructions for cospans in a quasi-category with finite colimits, denoted by $\text{coSpan}(\cQ)$.
If $\sdot$ denotes Waldhausen's construction for stable quasi-categories (as in \cite{barwickq}, \cite{barwickKtheory}, \cite{BGT}, or \cite{LurieHA}), Barwick and Rognes show in \cite[Proposition 3.7]{BarwickRognes} that these constructions are compatible with edgewise subdivision, in that there are levelwise equivalences of simplicial spaces
$$\text{Span}(\cQ)\simeq\esd(\sdot\cQ)\simeq \text{coSpan}(\cQ).$$
This gives an alternative proof that $\sdot\cQ$ is 2-Segal, and hence recovers the result from \cite{DK}.
\end{ex}

Unlike the flavors of $\sdot$-contructions discussed so far, the Waldhausen construction of an arbitrary Waldhausen category need not be a $2$-Segal space.  Let us look at a specific example.

\begin{ex}
\label{ex: Waldhausen categories}
Our criterion can be used to show that the $\sdot$-construction of the Waldhausen category of retractive spaces over $X$ need not be a $2$-Segal space, since its edgewise subdivision is not always a Segal space.  Let us look more closely at what goes wrong in this situation.

Given a space $X$, its {\em category of retractive spaces} is the category $R_f(X)$ whose objects are retractive spaces $(Z, i\colon X\to Z, r\colon Z\to X)$ over $X$, subject to the relative finiteness condition that $(Z, i(X))$ is a finite relative CW-complex, and maps that preserve the structure. As proven in \cite{waldhausen}, the category $R_f(X)$ can be given the structure of a Waldhausen category where cofibrations and weak equivalences are created in the
underlying category of CW-complexes. In the special case when $X$ is a singleton, we get the category $R_f(*)$ of retractive spaces over a point, which is just the category of finite pointed CW-complexes and pointed cellular maps, with cellular embeddings as cofibrations and homotopy equivalences as weak equivalences. 

In the context of the $\sdot$-construction for Waldhausen categories from \cite{waldhausen}, one can conclude from \cref{mainresult} that $\sdot R_f(X)$ is a $2$-Segal space if and only if its edgewise subdivision $Y(X):=\esd(\sdot R_f(X))$ is a Segal space.  We claim that if $X$ is a finite CW-complex, then the simplicial space $Y(X)$
is never a Segal space.  Here, we demonstrate that one of the Segal maps for $Y(*)$ is not an equivalence; the argument can be adapted for more general spaces $X$.

Let $\funsp$ be a finite $2$-dimensional CW-complex which is not contractible but whose suspension is contractible (for example, the classifying space of the perfect group from \cite[Example 2.38]{hatcher}). 
Using the notation $CP$ for the cone on $P$, consider the following diagrams in $R_f(*)$, which we denote by $D$ and $D'$:
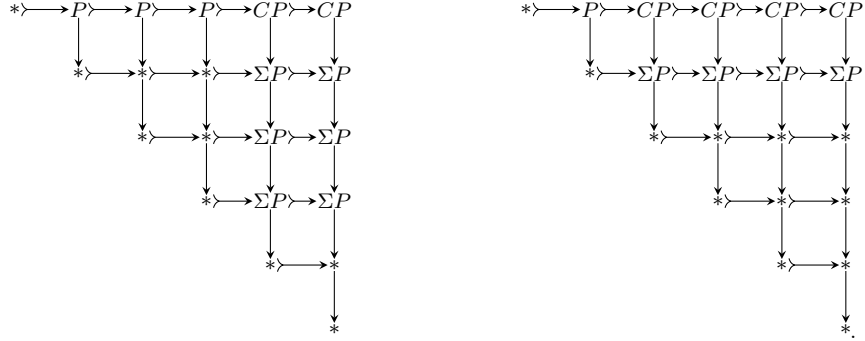
\[
\begin{tikzpicture}[scale=0.85, font=\footnotesize]%add , font=\scriptsize or , font=\small if necessary
\begin{scope}[outer sep=0, inner sep=0.4]
     \draw (1,0) node(a00){$*$};
\draw  (2,0)   node(a01){$\funsp$};
\draw (2,-1) node(a11) {$*$};
\draw  (3, -1)   node(a12){$*$};
\draw  (3, 0)   node(a02){$\funsp$};
\draw (3,-2) node (a22){$*$};
\draw  (4, 0)   node (a03){$\funsp$};
\draw  (4, -1)   node (a13){$*$};
\draw (4, -2)   node (a23){$*$};
\draw (4, -3) node (a33){$*$};
\draw  (5, 0)   node (a04){$C\funsp$};
\draw  (5, -1)  node (a14){$\Sigma\funsp$};
\draw  (5, -2)  node (a24){$\Sigma\funsp$};
\draw  (5, -3)  node (a34){$\Sigma\funsp$};
\draw (5, -4) node (a44){$*$};
\draw (6,0) node (a05){$C\funsp$};
\draw (6,-1) node (a15){$\Sigma\funsp$};
\draw (6,-2) node (a25){$\Sigma\funsp$};
\draw (6,-3) node (a35){$\Sigma\funsp$};
\draw (6,-4) node (a45){$*$};
\draw (6,-5) node (a55){$*$};

%horizontal
\draw[mono] (a00)--(a01);
\draw[mono] (a11)--(a12);
\draw[mono] (a01)--(a02);
\draw[mono] (a02)--(a03);
\draw[mono] (a12)--(a13);
\draw[mono] (a22)--(a23);
\draw[mono] (a03)--(a04);
\draw[mono] (a13)--(a14);
\draw[mono] (a23)--(a24);
\draw[mono] (a33)--(a34);
\draw[mono] (a04)--(a05);
\draw[mono] (a14)--(a15);
\draw[mono] (a24)--(a25);
\draw[mono] (a34)--(a35);
\draw[mono] (a44)--(a45);
%vertical
\draw[-stealth]  (a01)--(a11);
\draw[-stealth]  (a12)--(a22);
\draw[-stealth]  (a02)--(a12);
\draw[-stealth]  (a03)--(a13);
\draw[-stealth]  (a13)--(a23);
\draw[-stealth]  (a23)--(a33);
\draw[-stealth]  (a04)--(a14);
\draw[-stealth]  (a14)--(a24);
\draw[-stealth]  (a24)--(a34);
\draw[-stealth]  (a34)--(a44);
\draw[-stealth]  (a05)--(a15);
\draw[-stealth]  (a15)--(a25);
\draw[-stealth]  (a25)--(a35);
\draw[-stealth]  (a35)--(a45);
\draw[-stealth] (a45)--(a55);
\end{scope}

\begin{scope}[outer sep=0, inner sep=0.4, xshift=8cm]
     \draw (1,0) node(a00){$*$};
\draw  (2,0)   node(a01){$\funsp$};
\draw (2,-1) node(a11) {$*$};
\draw  (3, -1)   node(a12){$\Sigma\funsp$};
\draw  (3, 0)   node(a02){$C\funsp$};
\draw (3,-2) node (a22){$*$};
\draw  (4, 0)   node (a03){$C\funsp$};
\draw  (4, -1)   node (a13){$\Sigma\funsp$};
\draw (4, -2)   node (a23){$*$};
\draw (4, -3) node (a33){$*$};
\draw  (5, 0)   node (a04){$C\funsp$};
\draw  (5, -1)  node (a14){$\Sigma\funsp$};
\draw  (5, -2)  node (a24){$*$};
\draw  (5, -3)  node (a34){$*$};
\draw (5, -4) node (a44){$*$};
\draw (6,0) node (a05){$C\funsp$};
\draw (6,-1) node (a15){$\Sigma\funsp$};
\draw (6,-2) node (a25){$*$};
\draw (6,-3) node (a35){$*$};
\draw (6,-4) node (a45){$*$};
\draw (6,-5) node (a55){$*$};

%horizontal
\draw[mono] (a00)--(a01);
\draw[mono] (a11)--(a12);
\draw[mono] (a01)--(a02);
\draw[mono] (a02)--(a03);
\draw[mono] (a12)--(a13);
\draw[mono] (a22)--(a23);
\draw[mono] (a03)--(a04);
\draw[mono] (a13)--(a14);
\draw[mono] (a23)--(a24);
\draw[mono] (a33)--(a34);
\draw[mono] (a04)--(a05);
\draw[mono] (a14)--(a15);
\draw[mono] (a24)--(a25);
\draw[mono] (a34)--(a35);
\draw[mono] (a44)--(a45);
%vertical
\draw[-stealth]  (a01)--(a11);
\draw[-stealth]  (a12)--(a22);
\draw[-stealth]  (a02)--(a12);
\draw[-stealth]  (a03)--(a13);
\draw[-stealth]  (a13)--(a23);
\draw[-stealth]  (a23)--(a33);
\draw[-stealth]  (a04)--(a14);
\draw[-stealth]  (a14)--(a24);
\draw[-stealth]  (a24)--(a34);
\draw[-stealth]  (a34)--(a44);
\draw[-stealth]  (a05)--(a15);
\draw[-stealth]  (a15)--(a25);
\draw[-stealth]  (a25)--(a35);
\draw[-stealth]  (a35)--(a45);
\draw[-stealth] (a45)--(a55);

\draw (a55) node[xshift=0.1cm, yshift=-0.1cm]{$.$};
\end{scope}
\end{tikzpicture}
\] 
They define elements of $S_5(R_f(*))=\esd(S_2(R_f(*)))=Y(*)_2$. Consider their classes in $\pi_0(Y(*)_2)$.
We observe that, since $\funsp$ is not contractible, $D$ and $D'$ cannot define the same class in $\pi_0(Y(*))_2$, i.e.,
$$[D]\neq[D']\in\pi_0(Y(*))_2.$$
Their images under the $\pi_0$ of the Segal map of $Y(*)$,
$$\pi_0\left(\beta^2_1\right)\colon\pi_0(Y(*)_2)\to\pi_0\left(
Y(*)_1\overset{h}{\underset{Y(*)_0}{\times}}Y(*)_1
\right),$$
can be computed as the classes of the diagrams
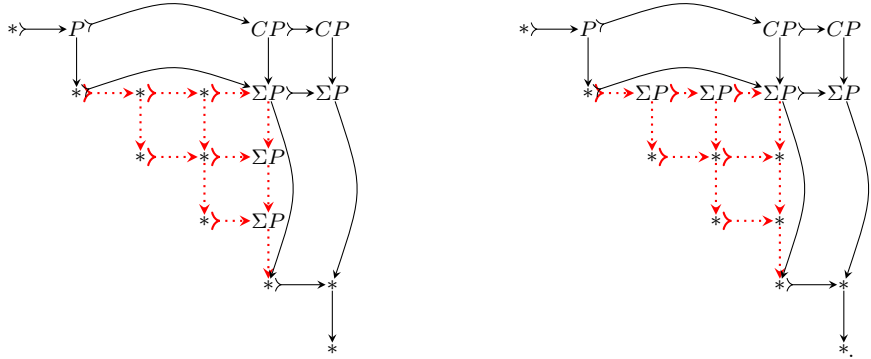
\tikzset{secondcomp/.style={dotted, thick, red}}
\[
\begin{tikzpicture}[scale=0.85, font=\footnotesize]%add , font=\scriptsize or , font=\small if necessary
\begin{scope}[outer sep=0, inner sep=0.4]
     \draw (1,0) node(a00){$*$};
\draw  (2,0)   node(a01){$\funsp$};
\draw (2,-1) node(a11) {$*$};
\draw  (3, -1)   node(a12){$*$};
%\draw  (3, 0)   node(a02){$\funsp$};
\draw (3,-2) node (a22){$*$};
%\draw  (4, 0)   node (a03){$\funsp$};
\draw  (4, -1)   node (a13){$*$};
\draw (4, -2)   node (a23){$*$};
\draw (4, -3) node (a33){$*$};
\draw  (5, 0)   node (a04){$C\funsp$};
\draw  (5, -1)  node (a14){$\Sigma\funsp$};
\draw  (5, -2)  node (a24){$\Sigma\funsp$};
\draw  (5, -3)  node (a34){$\Sigma\funsp$};
\draw (5, -4) node (a44){$*$};
\draw (6,0) node (a05){$C\funsp$};
\draw (6,-1) node (a15){$\Sigma\funsp$};
%\draw (6,-2) node (a25){$\Sigma\funsp$};
%\draw (6,-3) node (a35){$\Sigma\funsp$};
\draw (6,-4) node (a45){$*$};
\draw (6,-5) node (a55){$*$};

%horizontal
\draw[mono] (a00)--(a01);
\draw[mono, secondcomp] (a11)--(a12);
%\draw[mono] (a01)--(a02);
%\draw[mono] (a02)--(a03);
\draw[mono, secondcomp] (a12)--(a13);
\draw[mono, secondcomp] (a22)--(a23);
%\draw[mono] (a03)--(a04);
\draw[mono, secondcomp] (a13)--(a14);
\draw[mono, secondcomp] (a23)--(a24);
\draw[mono, secondcomp] (a33)--(a34);
\draw[mono] (a04)--(a05);
\draw[mono] (a14)--(a15);
%\draw[mono] (a24)--(a25);
%\draw[mono] (a34)--(a35);
\draw[mono] (a44)--(a45);
%vertical
\draw[-stealth]  (a01)--(a11);
\draw[-stealth, secondcomp]  (a12)--(a22);
%\draw[-stealth]  (a02)--(a12);%
%\draw[-stealth]  (a03)--(a13);
\draw[-stealth, secondcomp]  (a13)--(a23);
\draw[-stealth, secondcomp]  (a23)--(a33);
\draw[-stealth]  (a04)--(a14);
\draw[-stealth, secondcomp]  (a14)--(a24);
\draw[-stealth, secondcomp]  (a24)--(a34);
\draw[-stealth, secondcomp]  (a34)--(a44);
\draw[-stealth]  (a05)--(a15);
%\draw[-stealth]  (a15)--(a25);
%\draw[-stealth]  (a25)--(a35);
%\draw[-stealth]  (a35)--(a45);
\draw[-stealth] (a45)--(a55);

%%%%%%%%%%%%%%%%%%%%%%%%%
%compositions
\draw[mono] (a01)..controls (3.5, 0.5)..(a04);
\draw[mono] (a11)..controls (3.5, -0.5)..(a14);
\draw[-stealth] (a14)..controls (5.5, -2.5)..(a44);
\draw[-stealth] (a15)..controls (6.5, -2.5)..(a45);

%\draw (a55) node[xshift=0.1cm, yshift=-0.1cm]{$.$};
\end{scope}
%
%%%%%%%%%%%%%%%%%%%%%%%%%%%%%%%%%%%%%%%%%%%%%%
%
%
\begin{scope}[outer sep=0, inner sep=0.4, xshift=8cm]
    \draw (1,0) node(a00){$*$};
\draw  (2,0)   node(a01){$\funsp$};
\draw (2,-1) node(a11) {$*$};
\draw  (3, -1)   node(a12){$\Sigma\funsp$};
%\draw  (3, 0)   node(a02){$C\funsp$};
\draw (3,-2) node (a22){$*$};
%\draw  (4, 0)   node (a03){$C\funsp$};
\draw  (4, -1)   node (a13){$\Sigma\funsp$};
\draw (4, -2)   node (a23){$*$};
\draw (4, -3) node (a33){$*$};
\draw  (5, 0)   node (a04){$C\funsp$};
\draw  (5, -1)  node (a14){$\Sigma\funsp$};
\draw  (5, -2)  node (a24){$*$};
\draw  (5, -3)  node (a34){$*$};
\draw (5, -4) node (a44){$*$};
\draw (6,0) node (a05){$C\funsp$};
\draw (6,-1) node (a15){$\Sigma\funsp$};
%\draw (6,-2) node (a25){$*$};
%\draw (6,-3) node (a35){$*$};
\draw (6,-4) node (a45){$*$};
\draw (6,-5) node (a55){$*$};

%horizontal
\draw[mono] (a00)--(a01);
\draw[mono, secondcomp] (a11)--(a12);
%\draw[mono] (a01)--(a02);
%\draw[mono] (a02)--(a03);
\draw[mono, secondcomp] (a12)--(a13);
\draw[mono, secondcomp] (a22)--(a23);
%\draw[mono] (a03)--(a04);
\draw[mono, secondcomp] (a13)--(a14);
\draw[mono, secondcomp] (a23)--(a24);
\draw[mono, secondcomp] (a33)--(a34);
\draw[mono] (a04)--(a05);
\draw[mono] (a14)--(a15);
%\draw[mono] (a24)--(a25);
%\draw[mono] (a34)--(a35);
\draw[mono] (a44)--(a45);
%vertical
\draw[-stealth]  (a01)--(a11);
\draw[-stealth, secondcomp]  (a12)--(a22);
%\draw[-stealth]  (a02)--(a12);
%\draw[-stealth]  (a03), secondcomp--(a13);
\draw[-stealth, secondcomp]  (a13)--(a23);
\draw[-stealth, secondcomp]  (a23)--(a33);
\draw[-stealth]  (a04)--(a14);
\draw[-stealth, secondcomp]  (a14)--(a24);
\draw[-stealth, secondcomp]  (a24)--(a34);
\draw[-stealth, secondcomp]  (a34)--(a44);
\draw[-stealth]  (a05)--(a15);
%\draw[-stealth]  (a15)--(a25);
%\draw[-stealth]  (a25)--(a35);
%\draw[-stealth]  (a35)--(a45);
\draw[-stealth] (a45)--(a55);
%%%%%%%%%%%%%%%%%%%%%%%%%
%compositions
\draw[mono] (a01)..controls (3.5, 0.5)..(a04);
\draw[mono] (a11)..controls (3.5, -0.5)..(a14);
\draw[-stealth] (a14)..controls (5.5, -2.5)..(a44);
\draw[-stealth] (a15)..controls (6.5, -2.5)..(a45);

\draw (a55) node[xshift=0.1cm, yshift=-0.1cm]{$.$};
\end{scope}

\end{tikzpicture}
\]
Since the suspension of $\funsp$ is contractible, the diagrams $\beta^2_1(D)$ and $\beta^2_1(D')$ 
define the same class as the class of the constant diagram at a point, i.e.,
$$\pi_0(\beta^2_1)([D])=[\beta^2_1(D)]=[\beta^2_1(D')]=\pi_0(\beta^2_1)([D'])\in\pi_0\left(Y(*)_1\overset{h}{\underset{Y(*)_0}{\times}}Y(*)_1\right).$$
As a consequence, the function $\pi_0(\beta^2_1)$ is not injective, so the Segal map
$$\beta^2_1\colon Y(*)_2\to Y(*)_1\overset{h}{\underset{Y(*)_0}{\times}}Y(*)_1$$
cannot be an equivalence, and therefore $Y(*)$ cannot be a Segal space.
\end{ex}

\bibliographystyle{alpha}
\bibliography{ref}

\end{document}